\documentclass[twoside,draft,reqno]{birkart}

\usepackage{amssymb}

\input xy
\xyoption{all}

\begin{document}

\setcounter{page}{1}

%
%

\newcommand{\bo}[1]{\mathbf{#1}}
\newcommand{\lra}{\longrightarrow}
\newcommand{\ra}{\rightarrow}
\newcommand{\spaces}{\bo{Spaces}}
\newcommand{\algch}{\bo{Alg}_h^{\C}}
\def\alh#1{\bo{Alg}_h^{#1}}
\def\alg#1{{\bo{Alg}}^{#1}}
\newcommand{\algc}{\bo{Alg^C}}
\newcommand{\algt}{\bo{Alg^T}}
\newcommand{\T}{{\bo T}}
\newcommand{\C}{{\bo C}}
\newcommand{\X}{{\bo X}}
\newcommand{\A}{{\bo A}}
\newcommand{\D}{{\bo D}}
\newcommand{\bD}{\bar{\D}}
\newcommand{\bdn}{\bar{D}_n}
\newcommand{\Prj}{{\bo P}}
\newcommand{\M}{{\bo M}}
\newcommand{\gr}{\bo {Gr}}
\newcommand{\CC}{{\bar{\C}}}
\newcommand{\FC}{{\bo F_\ast \C}}
\newcommand{\bFC}{\overline{\FC}}
\newcommand{\HoM}{{{\bo H}{\bo o}\M}}
\newcommand{\sC}{{\spaces^{\C}}}
\newcommand{\st}{\spaces^{\T}}
\newcommand{\spr}{\spaces^{\Prj}}
\newcommand{\Tot}{{\rm Tot\,}}
\def\spc#1{\spaces^{#1}_{cof}}
\def\spf#1{\spaces^{#1}_{fib}}
\newcommand{\stf}{\spaces^{\T}_{fib}}
\newcommand{\stc}{\spaces^{\T}_{cof}}
\newcommand{\scf}{\spaces^{\C}_{fib}}
\newcommand{\scc}{\spaces^{\C}_{cof}}
\newcommand{\sFC}{\spaces^{\FC}}
\newcommand{\sFCc}{\spaces^{\FC}_{cof}}
\newcommand{\sprc}{\spaces^{\Prj}_{cof}}
\newcommand{\sprf}{\spaces^{\Prj}_{fib}}
\newcommand{\Cdi}{\C^{\rm disc}}
\newcommand{\sCdi}{\spaces^{\Cdi}}
\newcommand{\lsc}{{\bo L\sC}}
\newcommand{\lst}{{{\bo L}\st}}
\def\mapright#1#2{\smash{\mathop{\hbox to
#1pt{\rightarrowfill}}\limits^{#2}}}
\newcommand{\Hom}{{\rm Hom}}
\newcommand{\Map}{{\rm RMap}}
\newcommand{\map}{{\rm Map}}
\newcommand{\semith}{\bo{SemiTh}}
\newcommand{\algth}{\bo{AlgTh}}
\def\colim#1{{\rm colim}_{#1}\,}
\newcommand\hocolim{{\rm hocolim\,}}
\newcommand\holim{{\rm holim\,}}

\newtheorem{theorem}{Theorem}[section]
\newtheorem{lemma}[theorem]{Lemma}
\newtheorem{corollary}[theorem]{Corollary}
\newtheorem{proposition}[theorem]{Proposition}
\newtheorem{claim}[theorem]{Claim}
\newtheorem{remark}[theorem]{Remark}
\newtheorem{example}[theorem]{Example}
\newtheorem{conjecture}[theorem]{Conjecture}
\newtheorem{definition}[theorem]{Definition}
\newtheorem{problem}[theorem]{Problem}
\newtheorem{notation}[theorem]{Notation}
\newtheorem{note}[theorem]{Note}
\newtheorem{nn}[theorem]{}

\newenvironment{Proof}{\removelastskip\par\medskip
\noindent{\em Proof.} \rm}{\penalty-20\null\hfill${\,\lower0.9pt
\vbox{\hrule \hbox{\vrule height 0.3 cm \hskip 0.3 cm \vrule height 0.3
cm}\hrule}\,} $\par\medbreak}

\newenvironment{Proofx}{\removelastskip\par\medskip
\noindent{\em Proof.} \rm}{\par}

\newenvironment{Proofa}{\removelastskip\par\medskip
\noindent{\em Proof} \rm}{\penalty-20\null\hfill${\,\lower0.9pt
\vbox{\hrule \hbox{\vrule height 0.2 cm \hskip 0.2 cm \vrule height 0.2
cm}\hrule}\,} $\par\medbreak}

\clubpenalty = 10000
\widowpenalty = 10000

%
%

%
%

\title[From $\Gamma$-spaces to algebraic theories]
{From \boldmath$\Gamma$-spaces to algebraic theories\\
}

\author[B. Badzioch]{Bernard Badzioch}

\address{
Department of Mathematics\br
University of Minnesota\br
127 Vincent Hall\br
206 Church St. SE\br
Minneapolis, MN 55455\br
USA}

\email{badzioch@math.umn.edu}

\maketitle

\date {Sept. 24, 2001}

%
%

\section{Introduction}
\label{INTRO}
In \cite{segal} Segal gave the following elegant characterization  of
infinite loop spaces. 
Let $\bo \Gamma^{op}$ be the category whose objects are finite 
sets $[n]=\{0,1,\dots,n\}$ for $n\geq 0$ and whose morphisms are
all maps of sets $\varphi\colon [n]\ra [m]$ satisfying 
$\varphi(0)=0$. For $n>1$ and $1\leq k \leq n$ let 
$p^n_k\colon [n]\ra [1]$ be the map such that $p^n_k(i)=1$ only 
if $i=k$. 
A $\Gamma$-space  
is a functor $X\colon{\bo\Gamma^{op}}\ra\spaces$ such that 
$X[0]$ is a contractible space and for $n>1$ the map 
$\prod_{k=1}^n X(p^n_k)\colon X[n]\lra X[1]^n$ is a weak 
equivalence. 
The main result of \cite{segal} states that giving 
a $\Gamma$-space amounts to describing an infinite loop spaces 
structure on the space $X[1]$.

The general approach of Segal's paper proved to be useful to 
characterize a variety of homotopy invariant structures on spaces. 
In preparation for his work on infinite loop spaces Segal himself
showed that it is suitable for describing $A_\infty$-spaces.
Subsequently it was used by Bousfield \cite{bousfield} to give a
characterization  of $n$-fold loop spaces and by the author
\cite{badziochII} to identify generalized Eilenberg--Mac Lane spaces. 
The formalism underlying these examples can be described  as follows. 

\begin{definition}
\label{DEF OF SEMI}
A pointed semi-theory $\C$ is a category with objects 
$[0],[1],\dots$
such that $[0]$ is both an initial and a terminal object 
and for every $n>1$ there is a distinguished ordered set of $n$ 
different morphism $p^n_1,\dots,p^n_n\in\Hom_\C([n],[1])$.

\vskip .2cm

\noindent A homotopy algebra over a pointed semi-theory $\C$ is a functor 
$$X\colon\C\lra\spaces$$
such that $X[0]$ is a contractible space and the map 
$$\prod_{k=1}^n X(p^n_k)\colon X[n]\lra X[1]^n$$
is a weak equivalence for all $n>1$.
\end{definition} 

\noindent We will call the morphisms $p^n_k$ projection morphisms. 
It will be also convenient to denote by $p^1_1$ the identity 
morphism on the object $[1]$.

Given this definition a prototypical theorem about homotopy 
algebras would say that if $X$ is a homotopy algebra over  
a specified semi-theory $\C$ then the space $X[1]$ is equipped with 
some structure, depending only on $\C$. All results mentioned 
above are variants of this statement. 

It will be useful for us to consider also an unpointed version 
of the definition \ref{DEF OF SEMI} where we forget about the 
object $[0]$:

\begin{definition}
An unpointed semi-theory is a category $\C$ with objects 
$[1],[2],\dots$ and for every $n\geq 1$ a choice of projection 
maps $p^n_1,\dots,p^n_n\in\Hom_\C([n],[1])$.

A homotopy algebra over $\C$ is a functor $X\colon\C\ra\spaces$
such that for every $n>1$ the map $\prod_{k=1}^n X(p^n_k)$
gives a weak equivalence of spaces $X[n]$ and $X[1]^n$. 
\end{definition}

The major problem one faces while working with homotopy algebras is,
given a specific structure on a space, how to decide which semi-theory 
(if any) can be used to describe it. Or, going in the 
opposite direction, given a semi-theory $\C$ how to identify 
the structure it puts on the space $X[1]$ of a homotopy 
$\C$-algebra X.  
In \cite{badziochI} we showed that such questions are easier to settle 
if one restricts attention to the following  special kind of semi-theories:

\begin{definition}
A pointed (or unpointed) algebraic theory $\T$ is a pointed 
(resp. unpointed) semi-theory such
that  for $n>1$ the projections in $\T$ induce isomorphisms
$$\prod_{k=1}^n p^n_k\colon \Hom_\T([m],[n])\lra\Hom_\T([m],[1])^n$$
\end{definition}

\begin{definition}
Let $\C$ be a semi-theory. 
A strict $\C$-algebra is a functor 
$$A\colon\C\lra\spaces$$
such that the maps 
$$\prod_{k=1}^n A(p^n_k)\colon A[n]\lra A[1]^n$$
are isomorphisms for all $n>1$.
If $\C$ is pointed then we also assume that $A[0] = \ast$.
\end{definition}

The main result of \cite{badziochI} states that if $\T$
is an algebraic theory then every homotopy $\T$-algebra 
can be replaced by a weakly equivalent strict $\T$-algebra. 
On the other hand giving a strict $\T$-algebra $A$ amounts 
to providing the space $A[1]$ with some sort of an algebraic
structure determined by $\T$, e.g. monoid, group, ring, 
Lie algebra ~{...} (see \cite{lawvere},\cite{borceuxII}).
Therefore homotopy algebras over 
algebraic theories essentially  describe spaces weakly 
equivalent to algebraic objects in $\spaces$. 

In the present paper we show that it is possible to replace 
homotopy algebras by strict algebras even if one works with 
a semi-theory which is not an algebraic theory (as it is the
case with Segal's $\Gamma^{op}$) as follows.

\begin{theorem}
\label{MAIN 1}
For any pointed (or unpointed) semi-theory $\C$ 
there exists a poin\-ted (resp. unpointed) algebraic
theory  $\bFC$ such that the homotopy category of homotopy 
$\C$-algebras is equivalent to the homotopy category 
of strict $\bFC$-algebras. Moreover, the construction 
of $\bFC$ is functorial in $\C$. 
\end{theorem}

As an application of theorem \ref{MAIN 1} 
we can formulate conditions when two semi-theories describe 
equivalent structures on spaces.

\begin{theorem}
\label{MAIN 2}
A functor of semi-theories 
$G\colon\C\lra\C'$ induces an equivalence of 
the homotopy categories of homotopy algebras iff 
the induced functor 
$G\colon\bFC\ra\overline{\bo F_\ast\C'}$ 
between the associated algebraic
theories is a weak equivalence of categories, 
that is, G is a weak equivalence on the simplicial sets of 
morphisms of these categories.  
\end{theorem}

\vskip .5cm

\noindent{\bf Relationship to the theorem of Dwyer and Kan.\ }
 Let $\C$ be a small category and let $\D$ be its subcategory. 
Denote by $\spaces^{\C,\D}$ the category of all  diagrams
$X\colon\C\ra\spaces$ such that for every  morphism $\varphi\in\D$ the map
$X(\varphi)$ is a weak equivalence. Let $\D'\subseteq\C'$ be another pair of
categories and assume that we have a functor $G\colon (\C,\D)\ra(\C',\D')$
which is an epimorphism on the sets of objects. 
In \cite{dwyer-kanI} Dwyer and Kan gave sufficient and 
necessary conditions for the functor $G$ which guarantee 
that the induced functor 
$G^\ast\colon\spaces^{\C',\D'}\ra\spaces^{\C,\D}$ yields 
an equivalence of the homotopy categories of diagrams. 
Their approach involves constructing for a pair 
$(\C, \D)$ a simplicial category ${\bo F_\ast \C[\bo F_\ast \D]}$.
Since this construction is functorial, given $G$ as above
we obtain a new functor
$$FG\colon {\bo F_\ast \C[\bo F_\ast \D]}
\lra{\bo F_\ast \C'[\bo F_\ast \D']}$$ 
Then $G^\ast$ yields an equivalence  
of the homotopy categories of diagrams if and only if
the functor $FG$ is a weak equivalence of categories.

We note that our construction of the category $\bFC$
is very analogous to the construction of 
${\bo F_\ast \C[\bo F_\ast \D]}$. This is not very surprising
since homotopy algebras resemble the diagrams 
considered by Dwyer and Kan: a homotopy $\C$-algebra $X$
is a diagram of spaces satisfying the condition 
that certain maps defined in terms of $\C$ are weak equivalences. 
The difference is that in our setting the maps 
which are weak equivalences are not images of 
morphisms of $\C$. From this persprctive one can view 
theorem \ref{MAIN 2} as an extension of \cite[2.5]{dwyer-kanI}.

\vskip .5cm
\noindent{\bf Organization of the paper.\ }
Throughout this paper we will work with 
unpointed versions of theorems \ref{MAIN 1} and \ref{MAIN 2}.
The main difference between the pointed and the unpointed 
case is that the notion of a free unpointed semi-theory  
(sec. ~\ref{COMPLETION}) is more natural and as a consequence
the construction of the algebraic theory $\bFC$ is 
easier to describe. In section \ref{POINTED CASE} we comment  
on the changes required in order to perform 
this construction in the pointed case.
Beyond that our arguments work for 
pointed semi-theories without any major changes.
Hence, from now on
by semi-theory we will understand an unpointed semi-theory
(and similarly for algebraic theories).
   
We use freely the language of model 
categories. Section \ref{MOD CAT} contains a concise presentation of 
the main model category structures we will need.
In section ~\ref{COMPLETION},  we describe a functorial construction 
which associates to every semi-theory  $\C$ an algebraic theory 
$\bar{\C}$ (which we call the completion of $\C$ ) in such a way, 
that the categories of strict algebras 
over $\C$ and $\bar{\C}$ and canonically isomorphic. This construction,
however, usually does not preserve the category of homotopy algebras. 
To deal with this problem in section ~\ref{RESOLUTION} we replace 
the semi-theory $\C$ with $\FC$ - the simplicial resolution of
$\C$. This resolution is a simplicial semi-theory, that is 
a simplicial object in the category of semi-theories. Applying
the completion of \S\ref{COMPLETION} to every simplicial 
dimension of $\FC$ we obtain a simplicial algebraic theory 
$\bFC$. Then the results of section \ref{MOD CAT} show that 
the proof of theorem \ref{MAIN 1} amounts to verifying that certain
mapping complexes of $\FC$--diagrams are weakly equivalent. 
Assuming this we give the proof of theorem \ref{MAIN 2}.
The following two sections contain a few technical facts needed
to complete the proof of theorem \ref{MAIN 1}. 
In section \ref{INITIAL} we describe some properties 
of homotopy algebras over a semi-theory $\Prj$ which is 
the initial object in the category of semi-theories.
Then, in \S\ref{MAPS}, we deal with the problem of describing 
the space of maps between diagrams indexed by a simplicial 
category in terms of maps of diagrams over  discrete categories.  
Finally in sections \ref{PF LEMMA} and \ref{MAIN STEP} we complete 
the proof of theorem \ref{MAIN 1}. 

\vskip .5cm
\noindent{\bf Notation. \ }

(i)\ This paper is written simplicially. 
Thus, $\spaces$ stands for the
category of simplicial sets and consequently 
by  'space' we always mean is a simplicial set. 

(ii)\ For a small category $\C$ by $\sC$ we will denote 
the category of all functors $\C\ra\spaces$. 
If $X,Y$ are objects of $\sC$ then by $\Hom_\C(X,Y)$
we will mean the set of all natural transformations 
$X\ra Y$. Since objects of the indexing category 
$\C$ will be $[0],[1],\dots$ (except for \S\ref{MAPS})
this notation should not lead to confusing 
$\Hom_\C(X,Y)$ with $\Hom_\C([n],[m])$ - the set of morphisms 
$[n]\ra [m]$ in $\C$. 

\vskip .5cm  

\noindent{\bf Acknowledgment. \ } I would like to thank 
W. G. Dwyer for many discussions which contributed to this work.

%
%

\section{Model categories}
\label{MOD CAT}

Below we describe model category structures which
we will use throughout this paper. Our setting 
is basically the same as in \cite{badziochI}, therefore we 
discuss it briefly and refer there for the details.

\vskip .5cm

{\noindent \bf Model category for homotopy algebras.\ }
Let $\C$ be a small category. 
We consider two model category structures on
$\sC$, denoted by $\scf$ and $\scc$ \cite[sec. 3]{badziochI}. 
Weak equivalences in both cases are objectwise
weak equivalences. Fibrations  in $\scf$ are objectwise 
fibrations and cofibrations in  $\scc$ are objectwise cofibrations. 
The third class of morphisms in each of these model categories is 
as usual determined by the above choices.  
Both $\scf$ and $\scc$ are simplicial model categories and they
share the same simplicial structure: for $X\in\sC$ and a simplicial 
set $K$ the functor $X\otimes K$ is defined by 
$X\otimes K(c) = X(c)\times K$ for all $c\in\C$. 

For $X,Y\in\sC$ by $\map_\C(X,Y)$ we will denote the simplicial 
function complex of $X$ and $Y$, that is the simplicial set 
whose $k$--dimensional simplices are all natural 
transformations $X\otimes \Delta[k]\ra Y$. 

If $\C$ is a semi-theory then similarly as in 
\cite[sec. 5]{badziochI} we get an additional 
model category structure on $\sC$ which 
we will denote by $\lsc$. Cofibrations 
in $\lsc$ are the same as cofibrations in $\scf$. A map 
$f\colon X\ra Y$ is a weak equivalence in $\lsc$ if for 
every homotopy $\C$-algebra $Z$ which is fibrant as an object 
of $\scc$ the induced map of simplicial mapping complexes
$$f^\ast\colon\map_\C(Y, Z)\lra\map_\C(X,Z)$$
is a weak equivalence. To distinguish such maps from 
objectwise weak equivalences we will call them local equivalences. 

The significance of the model category $\lsc$ for our purposes
lies in the following fact \cite[5.7]{badziochI}:
\begin{proposition} 
The homotopy category of $\lsc$ is equivalent to the category 
obtained by taking the full subcategory of $\sC$ spanned by  
homotopy $\C$-algebras and inverting all objectwise 
weak equivalences. 
\end{proposition}
\noindent Thus, $\lsc$ is a model category suitable for describing 
the homotopy theory of homotopy $\C$-algebras.

\vskip .5cm 

{\noindent \bf The category of strict algebras. \ }
Let $\C$ be again a semi-theory and let $\alg{\C}$ denote 
the full subcategory 
of $\sC$ whose objects are all strict $\C$--algebras.
We have  the inclusion functor 
$J_{\C}\colon\alg{\C}\lra\sC$. 
Arguments analogous to these of
\cite[2.4]{badziochI} imply 

\begin{proposition}
There exists a functor 
$$K_{\C}\colon\sC\lra\alg{\C}$$
left adjoint to $J_{\C}$.
\end{proposition}

The category $\alg{\C}$ can be given
a simplicial model category structure with objectwise 
weak equivalences and objectwise fibrations. Using this 
fact we get \cite[6.3]{badziochI}

\begin{proposition} 
\label{KC JC}
The adjoint pair of functors $(K_{\C},J_{\C})$ is a Quillen 
pair between model categories $\lsc$ and $\alg{\C}$.
\end{proposition} 

\vskip .5cm 

{\noindent \bf Comparison lemma. \ }
The main step in the proof of theorem \ref{MAIN 1}
is to demonstrate that for certain semi-theories $\C$ the Quillen 
pair of (\ref{KC JC}) is a Quillen equivalence 
(which implies that the homotopy theories 
of homotopy $\C$-algebras and strict $\C$-algebras are 
equivalent). The next lemma shows that this in turn  
amounts to verifying that certain maps in 
$\lsc$ are local equivalences. 
For a semi-theory $\C$ denote by $C_n\in\sC$ the functor 
corepresented by $[n]$; that is 
$$C_n[m] := \Hom_{\C}([n],[m])$$ 
Let $\eta_{C_n}\colon C_n\lra J_\C K_\C C_n$ be the unit of the
adjunction $(K_\C, J_\C)$.

\begin{lemma}
\label{COMPARISON LEMMA}
If the map $\eta_{C_n}$ is a local equivalence for all $n\geq 0$ 
then the Quillen pair $(K_\C, J_\C)$ is a Quillen equivalence 
of the model categories $\lsc$ and $\alg{\C}$. 
\end{lemma} 
 
\begin{proof}
The arguments used in the proof of \cite[6.4]{badziochI}
show that $(K_\C, J_\C)$ is a Quillen equivalence if 
$\eta_X\colon X\lra J_\C K_\C X$ is a local equivalence
for any cofibrant object $X\in \lsc$. By assumption 
the map $\eta_X$ is a local equivalence if $X=C_n$.
We claim that $\eta_X$ is a local equivalence also if
$X$ is a finite disjoint sum of such functors,
$X= \coprod_{i}^{m} C_{n_i}$.
Indeed, notice that for $m\geq 1$ we have a map
$$\coprod_{j=1}^{m} C_1\lra C_m$$
which is induced by the projections $p^m_j\in\C$.
We can use it to construct a commutative diagram
$$\xymatrix{
\coprod_{i=1}^m C_{n_i}\ar[r]^\eta & 
J_\C K_\C(\coprod_{i=1}^n C_{n_i})\\
\coprod_{i=1}^m \coprod_{j=1}^{n_i} C_1\ar[u]^f\ar[d]_g\ar[r]^\eta &
J_\C K_\C(\coprod_{i=1}^m \coprod_{j=1}^{n_i} C_1)
\ar[u]_{J_\C K_\C(f)}\ar[d]^{J_\C K_\C(g)} \\
C_{\sum{n_i}}\ar[r]^\eta & J_\C K_\C(C_{\sum{n_i}}) 
}$$  
Directly from the definition of a local equivalence it follows 
that both left vertical maps are local equivalences.
One can also check that if $A$ is a strict $\C$-algebra then 
the maps induced by $f$ and $g$ on sets $\Hom_{\C}(-, A)$ 
are isomorphisms. This implies that the maps $J_\C K_\C(f)$ and 
$J_\C K_\C(g)$ are isomorphisms. 
By our assumption the bottom 
map $\eta$ is a local equivalence. Thus, 
also the top map must be a local equivalence.
This proves our claim.   
    
In order to show that $\eta_X$ is a local equivalence 
for an arbitrary cofibrant object $X$ one can proceed from here
using the same arguments as in the proof of  \cite[6.5]{badziochI}.
\end{proof}

%
%

\section{Completion of semi-theories}
\label{COMPLETION}

\begin{definition}
\label{DEF FREE SEMI}
A free semi-theory is a semi-theory $\C$ which is free as a category, and 
whose  projection morphisms are among free generators of $\C$.
\end{definition}

\noindent Our goal in this section is to describe 
a functorial construction 
which associates to every free semi-theory $\C$ an algebraic 
theory $\CC$ such that the categories of strict algebras over 
$\C$ and over $\CC$ are isomorphic. 

The category $\CC$ is defined as follows. The objects of $\CC$ are 
the same as the objects of $\C$:  
${\rm{ob}}\CC = {\rm{ob}}\C =\{[1],[2],\dots\}$.
For $n\geq 1$ the set of morphisms 
$\Hom_\CC([n],[1])$ consists 
of directed trees $T$:
{\objectmargin={-1pt}
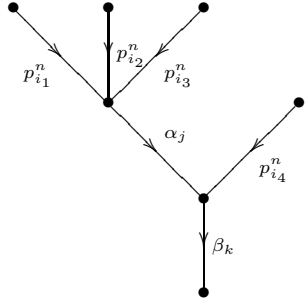
$$\xygraph{
!~:{@{-}|@{>}}
\bullet="v1"(:[dr]\bullet="v2"_{p^n_{i_1}} 
:[dr]\bullet="v3"^{\alpha_j}
:[d]\bullet="v4"^{\beta_k})&
\bullet(:"v2"^{p^n_{i_2}}) &
\bullet:"v2"^{p^n_{i_3}} & & \bullet:"v3"^{p^n_{i_4}}
}$$}
satisfying the following conditions:
\begin{itemize}
\item[1)] the lowest vertex of $T$ has only one incoming edge;
\item[2)] all edges of $T$ are labeled with $\alpha_i$  
where $\alpha\colon[m]\lra[k]$ is some generator 
of $\C$ and $1\leq i\leq k$
(if $\alpha = p^n_k$ then we will write $p^n_k$ rather than 
$(p^n_k)_1$); 
\item[3)] if a vertex of $T$ has $m$ incoming edges then the outgoing 
edge is labeled with $\alpha_i$ for some $\alpha\colon[m]\lra[k]$;
\item[4)] all the initial edges (that is, the edges starting at vertices 
with no incoming edges) are labeled with projections 
$p^n_k$ and no other edge of $T$ is labeled with 
any projection morphism.
\end{itemize}
For $m> 1$ we define
$$\Hom_\CC([n],[m]) := \Hom_\CC([n],[1])^m$$

If $(T_1,\dots,T_m)\in\Hom_\CC([n],[m])$ and $S\in\Hom_\CC([m],[1])$
then the composition $S\circ (T_1,\dots,T_m)$ is a tree obtained 
by grafting the tree $T_i$ in place of each initial edge of $S$ 
labeled with $p^m_i$. In general, if $(T_1,\dots,T_m)\in\Hom_\CC([n],[m])$
and $(S_1,\dots,S_k)\in\Hom_\CC([m],[k])$ then 
$$(S_1,\dots,S_k)\circ(T_1,\dots,T_m) := 
(S_1\circ(T_1,\dots,T_m),\dots, S_k\circ(T_1,\dots,T_m))$$
Notice that the identity morphism 
${\rm id}_{[n]}\in\Hom_\CC([n],[n])$ is given by 
${\rm id}_{[n]} = ({\bo p}^n_1,\dots , {\bo p}^n_n)$
where ${\bo p}^n_k$ is the tree

{\objectmargin={-1pt}
$$\xygraph{
!~:{@{-}|@{>}}
\bullet:[d]\bullet^{p^n_k} 
}$$}
We give $\CC$ a semi-theory structure by choosing the morphisms
${\bo p}^n_k$ to be projections in $\CC$. One can check 
that $\CC$ is, in fact,  an algebraic theory. 

Next we define a functor 
$$\Phi_\C\colon\C\lra\CC$$ which is the identity on objects, and 
such that $\Phi_\C(p^n_k) = {\bo p}^n_k$. If $\alpha\colon[n]\lra[m]$
is a generator of $\C$ and $\alpha$ is not a projection then 
$\Phi_\C(\alpha) = (T_{\alpha_1},\dots,T_{\alpha_m})$ where $T_{\alpha_k}$
is the tree
{\objectmargin={-1pt}
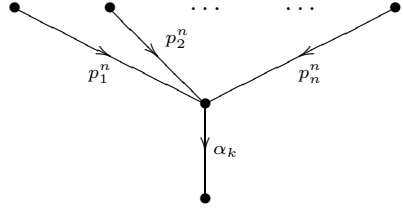
$$\xygraph{
!~:{@{-}|@{>}}
\bullet(:[drr]\bullet="u1"_{p^n_1} 
:[d]\bullet="w1"^{\alpha_k}) & 
\bullet(:"u1"^{p^n_2}) &
\dots & \dots &
\bullet(:"u1"^{p^n_n})
}$$}

\begin{definition}
\label{DEF COMPL}
We call the functor $\Phi_\C$ the completion of the semi-theory 
$\C$ to an algebraic theory.
\end{definition}

As we mentioned at the beginning of this section the essential property of 
the completion is that it preserves the category of strict algebras over 
$\C$:

\begin{proposition}
\label{COMPLETION ISO}
For any semi-theory $\C$ the functor $\Phi_\C\colon\C\lra\CC$
induces an isomorphism of the categories of strict algebras 
$$\Phi_\C^\ast\colon\alg{\CC}\mapright{20}{\cong}\alg{\C}$$
\end{proposition}

\noindent The proof will use the following
\begin{definition}
Let $T\in\Hom_\CC([n],[1])$. By $l(T)$ we will denote 
the number of edges of the longest (directed) path contained in 
$T$. If $T = (T_1,\dots,T_m)\in\Hom_\CC([n],[m])$ then 
$l(T) := \max\{l(T_1),\dots,l(T_m)\}$
\end{definition} 

\begin{Proofa}
{\em {\hskip -3mm} of proposition \ref{COMPLETION ISO}.\ }
We will construct a functor 
$$\Psi_\C\colon\alg{\C}\lra\alg{\CC}$$
which will be the inverse of $\Phi_\C^\ast$.
Let $X\colon\C\lra\spaces$ be a strict $\C$--algebra.
Set
$$\Psi_\C X[n]:=X[n]$$
For $n\geq 1$ denote by $\varphi_n\colon X[n]\lra X[1]^n$ 
the isomorphism induced by the product of projections,
$\varphi_n := \prod_k X(p^n_k)$. Let $T$ be a morphism of 
$\CC$. We define $\Psi_\C X(T)$ by induction with respect to 
$l(T)$:

\noindent 1) If $T\colon [n]\lra[1]$ and $l(T)= 1$ then 
$T= {\bo p}^n_k$ for some $1\leq k\leq n$. Then 
$$\Psi_\C X({\bo p}^n_k) := X(p^n_k)$$ 

\noindent 2) Assume that $\Psi_\C X(T)$ has been already defined 
for all $T\colon [n]\lra [1]$ such that $l(T)\leq r$ and let 
$(T_1,\dots, T_m)\in\Hom_\CC([n],[m])$, 
$l(T_1,\dots,T_m)\leq r$. Then for $x\in \Psi_\C X[n]$ 
set
$$\Psi_\C X(T_1,\dots, T_m )(x):= 
\varphi_m^{-1}(\Psi_\C X(T_1)(x),\dots, \Psi_\C X(T_m)(x))$$ 

\noindent 3) Let $\Psi_\C X(T)$ be defined for all morphisms
$T$ with $l(T)\leq r$ and let $S\in\Hom_\CC([n],[1])$, $l(S) = r+1$.
Then $S$ is of the form
{\objectmargin={0pt} 
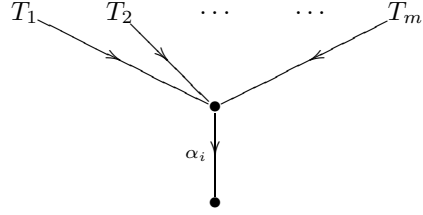
$$\xygraph{
!~:{@{-}|@{>}}
{T_1}(:[drr]\bullet="u1" 
:[d]\bullet="w1"_{\alpha_i}) & 
{T_2}(:"u1") &
\dots & \dots &
{T_m}(:"u1")
}$$}
where $(T_1,\dots,T_m)\colon [n]\lra [m]$ is a morphism with 
$l(T_1,\dots,T_m) = r$, $\alpha\colon [m]\lra[k]$ is a generator
of $\C$ and $1\leq i\leq k$. In this case 
$$\Psi_\C X(S) := X(p^m_i\circ\alpha)\circ\Psi_\C X(T_1,\dots, T_m)$$  
One can check that $\Psi_\C$ is a well defined functor and that 
the compositions $\Psi_\C\circ \Phi_\C^\ast$ and 
$\Phi_\C^\ast\circ\Psi_\C$ are both identities.
\end{Proofa}

\begin{remark}{\rm
While we will not use it, we note that completion of  a semi-theory 
can be described in categorical terms as follows. Let 
$\algth$ and $\semith$ denote the categories of algebraic theories 
and semi-theories respectively. We have the forgetful functor 
$$R\colon \algth\lra\semith$$
One can show that there exists a functor $L\colon\semith\lra\algth$
left adjoint to $R$. For a semi-theory $\C$ let 
$\eta_\C\colon\C\lra RL\C$ denote the unit of this adjunction. 
If $\C$ is a free semi-theory then we have $\CC = RL\C$
and $\Phi_C = \eta_\C$. Moreover, the functor induced by $\eta_\C$
on the categories of strict algebras is an isomorphism for any 
(not necessarily free) semi-theory $\C$.}
\end{remark}

%
%

\section{Simplicial resolution of a semi-theory}
\label{RESOLUTION}

Let $\C$ be a semi-theory. Following \cite[2.5]{dwyer-kanIII} 
we denote by $\FC$ the simplicial resolution of $\C$. 
Thus, $\FC$ is a simplicial category such that 
${\bo F_0\C}$ is a free category whose 
generators are morphisms of $\C$ and, for $k>0$, 
${\bo F_k\C}$ is a free category generated by $\bo F_{k-1}\C$.  
For every $[m],[n]\in\C$ there is a canonical map  $$\psi_{n,m}\colon
\Hom_\C([n],[m])\lra \Hom_{\FC}([n],[m])$$ We can define a semi-theory
structure on $\FC$ by choosing projection morphisms in 
$\FC$ to be the images of projections of $\C$ under the maps  
$\psi_{n,1}$.
We also have a functor \cite[2.5]{dwyer-kanIII}
$$\Psi\colon\FC\lra\C$$
Since for any $[m],[n]$ the composition
$$\Hom_\C([n],[m])\mapright{30}{\psi_{n,m}}
\Hom_{\FC}([n],[m])\mapright{30}{\Psi}\Hom_\C([n],[m])$$
is the identity map, the functor $\Psi$ is a map of semi-theories. 
Moreover the following holds:

\begin{proposition}
\label{C QE FC}
The functor $\Psi$ induces an adjoint pair of functors between 
model categories of homotopy algebras
$$\xymatrix{
\Psi_\ast\colon \lsc\ar@<0.5ex>[r]
&{\bo L\sFC}\colon\Psi^\ast\ar@<0.5ex>[l]
}$$
which is a Quillen equivalence.
\end{proposition}

\begin{proof}
By \cite[2.6]{dwyer-kanIII} the functor $\Psi$ is 
a weak equivalence of the categories $\C$ and $\FC$. 
It follows that the Quillen pair 
induced by $\Psi$ between the categories of diagrams
$$\xymatrix{ 
\Psi_\ast\colon\spf{\C}\ar@<0.5ex>[r]
&\spf{\FC}\colon\Psi^\ast\ar@<0.5ex>[l]
}$$
is a Quillen equivalence (this is essentially a consequence 
of \cite[2.1]{dwyer-kanIII}). 
The model categories $\lsc$
and $\bo L\sFC$ are obtained by localizing $\spf{\C}$ and 
$\spf{\FC}$, and so the statement follows from 
\cite[Thm. 3.4.20]{hirschhorn} which states 
that localization preserves Quillen equivalences. 

\end{proof}

Since for every $k\geq 0$ the category $\bo F_k\C$ 
is a free semi-theory  we can construct its completion 
to an algebraic theory (\ref{DEF COMPL}) 
$$ \Phi_k\colon {\bo F_k{\C}}\lra \overline{\bo F_k\C}$$
The functors $\Phi_k$ can be combined to define a functor 
of simplicial categories
$$\Phi\colon\FC\lra\bFC$$
where $\bFC$ is a simplicial algebraic theory which has 
$\overline{F_k\C}$ in its $k$-th simplicial dimension. 
Using (\ref{COMPLETION ISO}) we get

\begin{lemma}
\label{ALG ISO}
The functor $\Phi\colon\FC\lra\bFC$ 
induces an isomorphism of categories of strict algebras
$$\Phi^\ast\colon\alg{\bFC}\mapright{20}{\cong}\alg{\FC}$$
\end{lemma}

Recall (\ref{KC JC}) that we have a Quillen pair of functors
$$\xymatrix{
K_{\FC}\colon {\bo L\sFC}\ar@<0.5ex>[r]
&\alg{\FC}\colon J_{\FC}\ar@<0.5ex>[l]
}$$
In view of (\ref{C QE FC}) and (\ref{ALG ISO}) in order
to prove theorem \ref{MAIN 1} it is enough to show that the 
following holds.

\begin{proposition}
\label{QUILLEN PAIR}
The Quillen pair $(K_{\FC}, J_{\FC})$
is a Quillen equivalence. 
\end{proposition}

\noindent This last statement in turn is an immediate 
consequence of  lemma \ref{COMPARISON LEMMA}  
and 

\begin{lemma}
\label{FC UNIT}
For $n\geq 0$ let $FC_n\in\sFC$ denote the functor
corepresented by $[n]\in\FC$.
Then the unit of adjunction of the pair $(K_{\FC}, J_{\FC})$
$$\eta_{FC_n}\colon FC_n\lra J_{\FC}K_{\FC}FC_n$$
is a local equivalence.   
\end{lemma}

The proof of lemma \ref{FC UNIT} requires some technical preparations; 
we postpone it until \S\ref{PF LEMMA}. 
Meanwhile we turn to the proof
of theorem ~\ref{MAIN 2}. It will use the following fact
(see also \cite[8.6]{rezk}):

\begin{lemma}
\label{EQ FOR ALG}
Let $F\colon \T\ra\T'$ be a functor of algebraic theories. 
Then $F$ induces an adjoint pair of functors
between the categories of strict algebras 
$$\xymatrix{
F_\ast\colon \alg{\T}\ar@<0.5ex>[r] &
\alg{\T'} \colon F^\ast \ar@<0.5ex>[l] \\
}$$ 
which is a Quillen pair. Moreover, $F^\ast$ gives an 
equivalence of the homotopy theories of strict algebras
iff the functor $F$ is a weak equivalence of categories.
\end{lemma}

\begin{proof}
For the existence of the adjoint pair $(F_\ast, F^\ast)$
see \cite[3.7.7]{borceuxII}. It is a Quillen pair 
by \cite[3.4]{schwede}. Also by \cite[3.4]{schwede} 
if $F$ is a weak equivalence of $\T$ and $\T'$ 
then $(F_\ast,F^\ast)$ is a Quillen equivalence.
Therefore it remains to show that if $F^\ast$ induces 
an equivalence of the homotopy categories then 
$F$ is a weak equivalence. 

Let $T_n$ denote the $\T$-diagram corepresented 
by $[n]\in\T$. Since $\T$ is an algebraic theory, $T_n$ 
is a strict $\T$-algebra. Moreover, $T_n$ is a cofibrant
object in $\alg{\T}$. The functor $F^\ast$ gives 
an equivalence of the homotopy categories of strict
algebras, and thus its left adjoint $F_\ast$ provides 
the inverse equivalence. 
Therefore the unit of the adjunction $(F_\ast, F^\ast)$
$$\mu\colon T_n\lra F^\ast F_\ast T_n$$
is an (objectwise) weak equivalence in $\alg{\T}$. 
For $m\geq 0$
the map $\mu_{[m]}\colon T_n[m]\ra F^\ast F_\ast T_n[m]$
is given by maps of simplicial sets of morphisms 
$$T_n[m] = \Hom_\T([n],[m])\mapright{25}{F} \Hom_{\T'}([n],[m])
= F^\ast F_\ast T_n[m]$$
But this just means that 
$F\colon \T\ra\T'$ is a weak equivalence of categories.

\end{proof}

\begin{Proofa}
{\em {\hskip -3mm} of theorem \ref{MAIN 2}. \rm}
Let $G\colon \C\ra\C'$ be a functor of semi theories. 
Consider the diagram 
$$\xymatrix{
\alg{\bFC} \ar@<0.5ex>[r] \ar@<0.5ex>[d] &
\alg{\bFC'}\ar@<0.5ex>[l]\ar@<0.5ex>[d]\\
{\bo L\sFC} \ar@<0.5ex>[r]\ar@<0.5ex>[d]\ar@<0.5ex>[u]& 
{\bo L\spaces^{\FC'}}\ar@<0.5ex>[l]\ar@<0.5ex>[d]\ar@<0.5ex>[u]\\ 
\lsc \ar@<0.5ex>[r]\ar@<0.5ex>[u] & 
{\bo L\spaces^{\C'}} \ar@<0.5ex>[l]\ar@<0.5ex>[u]\\ 
}$$
in which every pair of arrows represents a Quillen pairs of functors. 
The horizontal pairs are induced by the functor $G$ while 
the vertical ones come from the adjunctions of (\ref{C QE FC}), 
(\ref{KC JC}) and (\ref{ALG ISO}). 
Propositions \ref{C QE FC} and \ref{QUILLEN PAIR} imply
that the vertical pairs are in fact Quillen equivalences. 
Therefore $G$ induces an equivalence of the homotopy categories of
$\lsc$ and  ${\bo L\spaces^{\C'}}$ if and only if it induces 
an equivalence of the homotopy categories of strict algebras
$\alg{\bFC}$ and $\alg{\bFC'}$. Thus lemma \ref{EQ FOR ALG}
completes the proof.

\end{Proofa}

%
%

\section{The initial semi-theory}
\label{INITIAL}

Denote by $\Prj$ the semi-theory whose only non-identity 
morphisms are projections. This is the initial object in 
the category of semi-theories: for any semi-theory $\C$
there is a unique map of semi-theories $\Prj\ra\C$.

The following fact states that local equivalences in
$\spaces^\Prj$ are particularly easy to detect.  

\begin{proposition}
\label{LOC WE IN PRJ}

A map $f\colon X\ra Y$ in $\spaces^\Prj$ is a local equivalence
if and only if the restriction $f_1\colon X[1]\ra Y[1]$
is a weak equivalence of spaces.

\end{proposition}

In the proof it will be convenient to use a definition of 
a local equivalence which is different 
(but equivalent, see \cite[18.6.31]{hirschhorn}) to 
the one given in \S\ref{MOD CAT}. First, 
assume that $X, Y$ are cofibrant in $\sprf$. Then 
a map $f\colon X\ra Y$ is a local equivalence if
for every homotopy algebra $Z$ fibrant in $\sprf$
the induced map of homotopy function complexes 
$$f^\ast\colon\map_\Prj(Y,Z)\lra\map_\Prj(X,Z)$$
is a weak equivalence of simplicial sets. 
If $X$ and $Y$ are not cofibrant then the map $f$
is a local equivalence if the induced map 
$f'\colon X'\ra Y'$ between cofibrant replacements 
of $X$ and $Y$ is one.

\begin{Proofa}
{\em {\hskip -3mm} of lemma \ref{LOC WE IN PRJ}.\rm\ }
Without loss of generality we can assume that $X$ and $Y$
are cofibrant in $\sprf$. Observe that for any strict 
$\Prj$-algebra $Z$ we have an isomorphism 
$$\map_\Prj(X, Z)\simeq \map(X[1], Z[1])$$
where the space on the right hand side is a mapping complex 
of simplicial sets. Moreover, given any Kan complex $K$ one can 
construct a strict $\Prj$-algebra $Z_K$ which is fibrant 
in $\sprf$and such that $Z_K[1]= K$. 
It easily follows from here that if 
$f\colon X\ra Y$ is a local equivalence then 
$f_1\colon X[1]\ra Y[1]$ is a weak equivalence. 

To see that the second implication holds 
notice that for any homotopy algebra $Z\in\spaces^\Prj$ 
one can find a strict $\Prj$-algebra $Z'$ 
such that $Z'$ is fibrant in $\sprf$
and there is an objectwise weak equivalence $Z\ra Z'$. 
 
\end{Proofa} 

In section \ref{PF LEMMA} we will use 
proposition \ref{LOC WE IN PRJ} to detect local 
equivalences between diagrams over a semi-theory $\FC$ 
(\S \ref{RESOLUTION}). In order to achieve that we will need 

\begin{lemma}
\label{QP FOR PRJ FC}
The map of semi-theories $J\colon\Prj\ra\FC$ induces 
Quillen pair of functors
$$\xymatrix{
J_\ast\colon \spaces^{\Prj}_{cof}\ar@<0.5ex>[r]&
\ \spaces^{\FC}_{cof}\colon J^\ast \ar@<0.5ex>[l] 
}$$

\end{lemma}   

\noindent The proof will follow from  

\begin{lemma}
\label{QP FOR PRJ C}
Let $J_\C\colon\Prj\ra\C$ denote the inclusion of $\Prj$
into a free semi-theory $\C$. Then the adjoint pair of functors
$$\xymatrix{
{J_\C}_\ast\colon \spaces^{\Prj}_{cof}\ar@<0.5ex>[r]&
\ \spaces^{\C}_{cof}\colon {J_\C}^\ast \ar@<0.5ex>[l] 
}$$
is a Quillen pair

\end{lemma}

\begin{proof}
The functor ${J_\C}_\ast\colon\spaces^{\Prj}\ra\spaces^{\C}$
can be described as follows. If $Y\in\spaces^{\Prj}$ and 
$[n]\in\C$ then 
$$ {J_\C}_\ast Y[n] = Y[n] \sqcup \coprod_{(\varphi\colon [m]\ra [n]) 
\in C_n} Y[m]$$
where $C_n$ is the set of all morphisms 
$\varphi\colon [m]\ra [n]$ such that 
$\varphi = \alpha_k\circ\alpha_{k-1}\circ\dots\circ\alpha_1$ 
where $\alpha_1, \dots, \alpha_k $ are generators of $\C$ and 
$\alpha_1$ is not a projection.

From this description it is clear that ${J_\C}_\ast$ preserves 
objectwise cofibrations and weak equivalences. Therefore 
$({J_\C}_\ast, {J_\C}^\ast)$ is a Quillen pair.

\end{proof}

\begin{Proofa}
{\em {\hskip -3mm} of lemma \ref{QP FOR PRJ FC}. \rm\ }
For every $k= 0,1,\dots$ the functor $J$ defines a map of
semi-theories
$$J_k\colon \Prj\lra \bo F_k\C $$ 
Let $s\spaces$ denote the category of simplicial spaces taken with the 
model structure in which  cofibrations and weak equivalences 
are the objectwise cofibrations and weak equivalences. 
The functors 
$$(J_k)_\ast\colon \sprc\lra \spaces^{F_k\C}_{cof}$$
can be assembled into a functor  
$$J_\bullet\colon \sprc\lra s\sFCc$$
such that 
$$\xymatrix{
J_\bullet Y[n] = ((J_0)_\ast Y[n] & 
(J_1)_\ast Y[n]\ar@2{->}[l] & \dots)\ar@3{->}[l]  \\
}$$
Take the diagonal functor
$$ |-|\colon s\spaces\lra\spaces$$
which associates to every simplicial space $X_\bullet$ its diagonal 
$|X_\bullet|$. It induces a functor 
$$|-|\colon s\sFCc \lra \sFCc$$
One can check that  $J_\ast\colon \sprc\ra\sFCc$ is equal to 
the composition  
$$\sprc \mapright{20}{J_\bullet} s\sFCc\mapright{20}{|-|}\sFCc$$
Since $\bo F_k\C$ is a free semi-theory for all $k\geq 0$, 
lemma \ref{QP FOR PRJ C} implies that  $J_\bullet$ 
sends cofibrations and weak equivalences from $\sprc$ 
to objectwise weak equivalences and cofibrations  in $s\sFCc$. 
The diagonal functor, in turn sends such maps to  weak equivalences 
and cofibrations in $\sFCc$. Therefore the functor
$J_\ast\colon\spaces^{\Prj}_{cof}\ra\spaces^{\FC}_{cof}$  
preserves cofibrations and  weak equivalences and  the adjoint pair 
$(J_\ast, J^\ast)$ is a Quillen pair as claimed.

\end{Proofa}

%
%

\section{Mapping complexes of diagrams over a simplicial category}
\label{MAPS}
Let $\C$ be an arbitrary simplicial category.   
Denote by $\C_k$ the category in $k$-th simplicial dimension
of $\C$. Our goal in this section is to show that given $X,Y\in\sC$ 
one can describe the simplicial mapping complex $\map_\C(X,Y)$ by 
means of mapping complexes of diagrams over the discrete 
categories $\C_k$. 

Thus, assume that $X$ is a diagram of spaces over $\C$.
We consider two ways in which we can turn $X$ into a diagram over 
the category $\C_k$:
\begin{itemize}
\item[1)] For $c\in\C_k$ let $X_k(c)$ be the set of 
$k$--dimensional simplices of $X(c)$ (regarded as a discrete 
simplicial set). 
If $\varphi\colon c\ra d$ is a morphism of $\C_k$ then 
the functor $X$ gives a map $X_k(c)\ra X_k(d)$. Thus 
we obtain a diagram $X_k\in\spaces^{\C_k}$. 
One can check that the correspondence $X\mapsto X_k$ 
defines a functor $R_k\colon\sC\ra\spaces^{\C_k}$

\item[2)] Let $\varphi\colon c\ra d$ again be a morphism of $\C_k$.
The diagram $X$ provides a map 
$$X(\varphi)\colon X(c)\times \Delta[k]\lra X(d)$$
This can be used to define a functor
$L_kX\colon \C_k\lra\spaces$ such that 
$L_kX(c):=X(c)\times\Delta[k]$ and for $\varphi\colon c\ra d$
in $\C_k$ the map $L_k X(\varphi)$ is given by 
$$L_k X(\varphi)\colon X(c)\times\Delta[k]
\mapright{45}{X(\varphi)\times pr_2} X(d)\times \Delta[k]$$
where $pr_2$ is the projection on the second factor. 
Again, we observe that this construction is functorial in 
$X$ and yields a functor $L_k\colon\sC\ra\spaces^{\C_k}$. 
\end{itemize}

Let $\theta\colon k\ra l$ be a morphisms in 
the simplicial indexing category $\Delta^{op}$. 
Since $\C$ is a simplicial
category we have a functor $\theta\colon\C_k\ra\C_l$ 
which induces 
$$\theta^\ast\colon\spaces^{\C_l}\lra\spaces^{\C_k}$$
For any $X\in\sC$ there are  
obvious natural transformations of $\C_k$--diagrams:
$$f_\theta\colon X_k\lra \theta^\ast X_l$$
and 
$$g_\theta\colon \theta^\ast L_l X \lra L_k X$$
Now, for  $X,Y\in\sC$ consider the function complex 
$\map_{\C_k}(X_k, L_k Y)$. A morphism 
$\theta\in\Delta^{op}$ as above defines a map 
$$\map_{\C_l}(X_l,L_l Y)\lra 
\map_{\C_l}(\theta^\ast X_l, \theta^\ast L_l Y)$$
Composing it with the maps induced by  
$f_\theta$ and $g_\theta$ yields a map 
$$\theta^\ast\colon\map_{\C_l}(X_l, L_l X)
\lra\map_{\C_k}(X_k, L_k Y)$$
It is not hard to check that it defines 
a cosimplicial space $\map_{\C_\bullet}(X,LY)$
which has $\map_{\C_k}(X_k, L_k Y)$ in $k$-th
cosimplicial dimension.
We claim that the simplicial function complex 
$\map_{\C}(X,Y)$ can be recovered from this cosimplicial 
space. 
Recall \cite[VIII.1 p.390]{goerss-jardine} that 
to every cosimplicial space $Z_\bullet$ one can associate 
its total space ${\rm Tot}(Z_\bullet)$. 

\begin{proposition}
\label{TOT LEMMA}
For any simplicial category $\C$ and $X,Y\in\sC$
the simplicial sets $\map_\C(X,Y)$ and 
$\Tot\map_{\C_\bullet}(X,LY)$ are isomorphic. 
Moreover, the isomorphism is functorial in $X$ and $Y$. 
\end{proposition}

\begin{proof}
The total space $\Tot\map_{\C\bullet}(X,LY)$
is the equalizer of the diagram 
$$\xymatrix{
\displaystyle{\prod_{k\geq 0}} \map_{\C_k}(X_k\times\Delta[k],L_k Y)
\ar@<1.5ex>[r]^{\delta_1}\ar@<0.5ex>[r]_{\delta_2} &
\displaystyle{\prod_{\phi\colon k\ra l}} 
\map_{\C_k}(X_k\times\Delta[l], L_k Y)
}$$
where the second product is indexed by all morphisms 
$\phi\in\Delta^{op}$. 
We claim that there exists a map 
$$\psi\colon \map_\C(X,Y)\lra 
\prod_{k\geq 0}\map_{\C_k}(X_k\times\Delta[k],L_k Y)$$
such that $\delta_1\circ\psi = \delta_2\circ\psi$.
As a consequence there is a unique map 
$\bar{\psi}\colon\map_\C(X,Y)\ra\Tot\map_{\C_\bullet(X,LY)}$
which fits into the commutative diagram 
$$\xymatrix{
\map_\C(X,Y)\ar[d]^{\bar\psi}\ar[rd]^{\psi} & \\
\Tot\map_{\C_\bullet}(X,LY)\ar[r] &
\prod_{k\geq 0}\map_{\C_k}(X_k\times\Delta[k], L_k Y)
}$$  
The map $\psi$ is defined as follows.
For $k \geq 0$ let 
$\epsilon_k\colon X_k\times\Delta[k]\ra L_k X$
be a natural transformation given by   
$$\epsilon_k = i_k\times pr_2 \colon 
X_k(c)\times\Delta[k]\lra X(c)\times\Delta[k]$$
where $i_k$ is the inclusion of the set of $k$-dimensional simplices 
into $X(c)$ and $pr_2$ is the projection on the second factor.
Consider the induced map of simplicial sets
$$\epsilon_k^\ast\colon\map_{\C_k}(L_kX,L_kY)\lra
\map_{\C_k}(X_k\times\Delta[k], L_k Y)$$
Its composition with the map 
$\map_\C(X,Y)\ra\map_{\C_k}(L_k X,L_k Y)$ given by 
the functor $L_k$ yields 
$$\psi_k\colon \map_\C(X,Y)\lra
\map_{\C_k}(X_k\times\Delta[k], L_k Y)$$
We define $\psi := \prod_{k\geq 0} \psi_k$.

One can check that $\psi$ is an embedding 
of simplicial sets and thus so is $\bar{\psi}$. It is also 
not hard to see that $\bar{\psi}$ is an epimorphism. 
Therefore it gives an isomorphism of $\map_\C(X,Y)$
and $\Tot\map_{\C_\bullet}(X,LY)$. 
  
\end{proof}

%
%

\section{Proof of lemma \ref{FC UNIT}.}
\label{PF LEMMA}
Let $\D$ be a free semi-theory and $\Phi_{\D}\colon\D\lra\bar{\D}$
be the completion of $\D$ to an algebraic theory. Denote by 
$D_n$ the $\D$ -- diagram corepresented by $[n]$.
Using the functor $\Phi_{\D}$ we can also define a $\D$ -- diagram 
$\bar{D}_n$ such that $\bar{D}_n[m] = \Hom_{\bar{\D}}([n],[m])$,
and a map of $\D$ -- diagrams 
$$\Phi_{\D}^\ast\colon D_n\lra \bdn$$
Moreover, since  $\Phi_{\D}$ is an embedding of categories,
$D_n$ is a subdiagram of $\bdn$. 
We define a filtration of the diagram $\bdn$ by $\D$--diagrams
$$\bdn^0\subseteq\bdn^1\subseteq\dots\bdn$$
as follows. Set $\bdn^0 := D_n$. If $\bdn^i$ is defined for $i\leq k$
then  $\bdn^{k+1}$ is the smallest $\D$--subdiagram of $\bdn$
such that if $T_1,T_2,\dots,T_m$ are elements of $\bdn^k[1]
\subseteq \Hom_{\bar{\D}}([n],[1])$ then 
$(T_1,T_2,\dots,T_m)\in\Hom_{\bar{\D}}([n],[m])$ belongs to 
$\bdn^{k+1}[m]$. From the definition of $\bar{\D}$ (\S\ref{COMPLETION})
it follows that  $\colim{k}\bdn^k = \bdn$. 

Recall (\S \ref{INITIAL}) that by $\Prj$ we denoted the semi-theory 
which has projections as the only non-identity morphisms. 
The unique map $\Prj\ra\D$ induces a $\Prj$-diagram structure 
on $D_n$, $\bdn$ and $\bdn^k$. 
Similarly as above we define a filtration of 
$\bdn$ by $\Prj$--diagrams 
$$s\bdn^0\subseteq s\bdn^1\subseteq\dots\subseteq\bdn$$
where $s\bdn^0 = D_n$ and $s\bdn^{k+1}$ is the smallest 
$\Prj$--subdiagram of $\bdn$ such that if 
$T_1,T_2,\dots,T_m\in \bdn^k[1]$ then 
$(T_1,T_2,\dots,T_m)\in s\bdn^{k+1}[m]$. 
We have inclusions of $\Prj$--diagrams
$\bdn^k\subseteq s\bdn^{k+1}\subseteq \bdn^{k+1}$  
and $\colim{k}s\bdn = \bdn$. 

We claim that the filtrations 
$\{\bdn^k\}$ and $\{s\bdn^k\}$ have the following
property:
 
\begin{lemma}
\label{PULLBACK LEMMA}
For any $\D$--diagram of spaces $X\colon\D\lra\spaces$
and for $k\geq 0$ the square of simplicial mapping complexes
$$\xymatrix{
\map_{\D}(\bdn^k,X)\ar[d] & \map_{\D}(\bdn^{k+1}, X)\ar[l]\ar[d] \\
\map_{\Prj}(\bdn^k, X)     & \map_{\Prj}(s\bdn^{k+1}, X)\ar[l] \\
}$$ 
is a pull-back diagram.
\end{lemma}

\noindent The proof of this fact is given in section \ref{MAIN STEP}.

Now we can proceed with the proof of lemma \ref{FC UNIT}. 
Consider $\bo F_m\C$ --  the free semi-theory in the $m$-th 
simplicial dimension of $\FC$, and let $\overline{\bo F_m\C}$ 
be its completion to an algebraic theory. 
Setting $\D := \bo F_m\C$ above 
we see that the $\bo F_m\C$-diagram $\overline{F_mC}_n$, 
(where $\overline{F_mC}_n[r]=\Hom_{\overline{\bo F_m\C}}([n],[r])$) 
admits two filtrations:
by $\bo F_m\C$-diagrams
$$ F_mC_n = \overline{F_mC}_n^0\subseteq \overline{F_mC}_n^1
\subseteq\dots\subseteq \overline{F_mC}_n$$
and by $\Prj$--diagrams
$$ F_mC_n = s\overline{F_mC}_n^0\subseteq s\overline{F_mC}_n^1
\subseteq\dots\subseteq \overline{F_mC}_n$$
The first of these filtrations, 
combined for all $m$, yields a filtration of the diagram 
$\overline{FC}_n$ by $\FC$--diagrams 
$$FC_n = \overline{FC}_n^0\subseteq \overline{FC}_n^1
\subseteq\dots\subseteq \overline{FC}_n$$
Similarly, the filtrations of $\overline{F_mC}$ by $\Prj$--diagrams
$s\overline{F_mC}_n^k$ for $m \geq 0$ give a filtration of 
$\overline{FC}_n$ by $\Prj$-diagrams
$$FC_n = s\overline{FC}_n^0\subseteq s\overline{FC}_n^1
\subseteq\dots\subseteq \overline{FC}_n$$.

\begin{lemma}
\label{DIAG}
For $X\in\spaces^{\FC}$ consider the following diagrams of 
simplicial function complexes:
$$\xymatrix{
\map_{\FC}(\overline{FC}_n^k,X)\ar[d] & 
\map_{\FC}(\overline{FC}_n^{k+1}, X)\ar[l]_f\ar[d] \\
\map_{\Prj}(\overline{FC}_n^k, X)  & 
\map_{\Prj}(s\overline{FC}_n^{k+1}, X)\ar[l]_g \\
}$$ 
This is a pullback diagram for all $X$, and $k,n\geq 0$.
\end{lemma}
 
\begin{proof}
Take the diagrams of cosimplicial spaces (sec. \ref{MAPS})
$$\xymatrix{
\map_{\bo F_\bullet\C}(\overline{FC}_n^k,LX)\ar[d] & 
\map_{\bo F_\bullet\C}(\overline{FC}_n^{k+1}, LX)
\ar[l]_{f_\bullet}\ar[d] \\
\map_{\Prj}(\overline{FC}_n^k, LX)  & 
\map_{\Prj}(s\overline{FC}_n^{k+1}, LX)\ar[l]_{g_\bullet} \\
}$$ 
Since limits of cosimplicial spaces can be computed by 
taking limits in each cosimplicial dimension separately, 
lemma \ref{PULLBACK LEMMA} implies that this is a pullback 
diagram of cosimplicial spaces. Therefore out claim is 
a consequence of (\ref{TOT LEMMA}) and of 
the fact that the functor $\Tot$ commutes with limits.
\end{proof}

Consider the map $g$ in the statement of 
lemma \ref{DIAG}. Our next goal is  

\begin{lemma}
\label{LOWER MAP}
Let $X$ be a homotopy algebra fibrant in $\spaces_{cof}^{\FC}$.
Then for every $k\geq0$
the map 
$$g\colon \map_{\Prj}(s\overline{FC}_n^{k+1}, X)\lra
\map_{\Prj}(\overline{FC}_n^k, X)$$
is an acyclic fibration of simplicial sets.
\end{lemma}
 
\begin{proof}
The map $g$ is induced by an inclusion 
$\iota_k\colon\overline{FC}_n^k\hookrightarrow s\overline{FC}_n^{k+1}$. 
Since all  inclusions are cofibrations in $\spaces_{cof}^{\FC}$
we get that $g$ is a fibration. 
Thus, it remains to show that $g$ is a weak equivalence 
of simplicial sets. 

By lemma \ref{QP FOR PRJ FC} if $X$ a homotopy $\FC$-algebra 
fibrant in $\spaces_{cof}^{\FC}$ then it is also a 
$\Prj$-homotopy algebra which is fibrant in $\spaces_{cof}^\Prj$. 
Therefore it is enough to show that the map $\iota_k$
is a local equivalence in $\spaces^\Prj$. This, however, 
is a consequence of (\ref{LOC WE IN PRJ}) and the observation that 
$\iota_k$ restricts to an isomorphism of simplicial sets
$$\overline{FC}_n^k[1]\mapright{25}{\cong} 
s\overline{FC}_n^{k+1}[1]$$

\end{proof}

Next, consider the upper map $f$ in the diagram in lemma \ref{DIAG}. 
Combining (\ref{LOWER MAP}) and (\ref{DIAG}) we obtain

\begin{corollary}
\label{AC FIB}
For all $n,k\geq 0$ the map 
$$f\colon\map_{\FC}(\overline{FC}_n^{k+1}, X)\lra
\map_{\FC}(\overline{FC}_n^k,X)$$
is an acyclic fibration of simplicial sets.
\end{corollary}

We use this fact to finish the proof of lemma \ref{FC UNIT}. 
The map  $\eta_{FC_n}\colon FC_n\ra J_\C K_\C FC_n$ 
is given by the inclusion of $\FC$-diagrams 
$FC_n = \overline{FC}_n^0\hookrightarrow \overline{FC_n}$. 
Moreover, $\overline{FC}_n$ is represented by the colimit 
$\overline{FC}_n = \colim{k}\overline{FC}_n^k$.
We have a commutative diagram 
$$\xymatrix{
FC_n \ar[dr]_{\eta_{FC_n}}\ar[rr] & & 
\hocolim_k \overline{FC}_n^k \ar[dl] \\
 & \colim{k}\overline{FC}_n^k & \\
}$$
where the homotopy colimit is taken in the model 
category $\bo L\sFC$. Corollary \ref{AC FIB} implies 
that both maps $FC_n\ra \hocolim_{k}\overline{FC}_n^k$
and $\hocolim_k \overline{FC}_n^k\ra\colim{k}\overline{FC}_n^k$ are
local equivalences. Therefore the map $\eta_{FC_n}$ is also 
a local equivalence.

%
%

\section{The pullback lemma}
\label{MAIN STEP}

In order to prove lemma \ref{PULLBACK LEMMA} we need 
to make a few observations  about $\bdn$ and its filtrations. 
For the rest of this section  $\alpha_i$ will always denote 
a generator of the free semi-theory  $\D$. In particular, 
if $\varphi$ is a morphism of $\D$  then by
$\varphi =\alpha_k\circ\dots\circ\alpha_1$ 
we will understand the decomposition of $\varphi$ into 
generators of $\D$. 

\begin{lemma} 
\label{L1}
Let $\varphi = \alpha_k\circ\dots\circ\alpha_1$ be a morphism of
$\D$ such that $\alpha_1$ is not a projection morphism. 
If $T, T'\in\bdn$ satisfy $\varphi(T) = \varphi(T')$ then 
$T=T'$.
\end{lemma} 

\begin{proof}
Assume that $\varphi=\alpha_1:[m]\lra[r]$ is a generator of $\D$. 
If $T\in\bdn[m]$ then $T=(T_1,\dots,T_m)$ where 
$T_1,\dots,T_m$ are trees contained in $\Hom_{\bD}([n],[1])$.
Using the definition of composition of morphisms in $\bD$ 
we get that $\varphi(T)= (\varphi_1 T,\dots, \varphi_k T)$
where $\varphi_i T$ is a tree in $\bD$ of the form
{\objectmargin={0pt} 
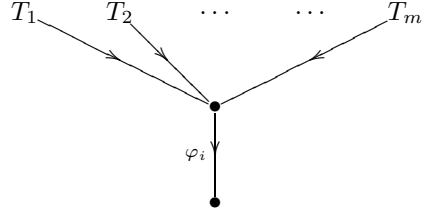
$$\xygraph{
!~:{@{-}|@{>}}
{T_1}(:[drr]\bullet="u1" 
:[d]\bullet="w1"_{\varphi_i}) & 
{T_2}(:"u1") &
\dots & \dots &
{T_m}(:"u1")
}$$}
Since $\varphi_i T$ contains all information about $T$
the equality $\varphi T$ = $\varphi T'$ must imply  
$T = T'$.

Next, assume that  $\varphi = \alpha_2\circ\alpha_1$ where 
$\alpha_1\colon[m]\lra[r]$ and $\alpha_2$ 
is a projection morphism, say $\alpha_2 = p^r_j$. In this case 
$\varphi(T)$ is the tree $(\alpha_1)_j T$: 
{\objectmargin={0pt} 
$$\xygraph{
!~:{@{-}|@{>}}
{T_1}(:[drr]\bullet="u1" 
:[d]\bullet="w1"_{(\alpha_i)_j}) & 
{T_2}(:"u1") &
\dots & \dots &
{T_m}(:"u1")
}$$}
Thus, as before, $\varphi T = \varphi T'$ implies that $T = T'$.

To show that the lemma holds for an arbitrary 
$\varphi=\alpha_k\circ\dots\circ\alpha_1$ one can now argue
by induction with respect to $k$.     
\end{proof}

The next fact is a relative version of (\ref{L1}).

\begin{lemma}
\label{L2}
Let $T, T'\in\bdn$ and $\varphi=\alpha_k\circ\dots\circ\alpha_1$,
$\varphi'=\alpha_{k'}'\circ\dots\circ\alpha_1'$ be morphisms 
of $\D$ such that $k\leq k'$ and $\alpha_1$, $\alpha_1'$ are not
projections. If $\varphi(T)= \varphi(T')$ then 
$\varphi'=\varphi\circ\theta$ and $T = \theta(T')$
where $\theta= \alpha_{k'-k}\circ\dots\circ\alpha_1$
\end{lemma}

\begin{proof}
We use induction with respect to $k$. 
If $k=0$ then $\varphi= {\rm id}$ and the claim obviously holds.
Assume that $k\geq 0$ and that $\alpha_k$  is not a projection. 
Directly  inspecting trees as in the proof of (\ref{L1})
one sees that the equality $\varphi(T)=\varphi(T')$
implies $\alpha_{k'}=\alpha_k$.
Therefore we have
$$\alpha_k\circ(\alpha_{k-1}\circ\dots\circ\alpha_1)(T)
=\varphi(T)=\varphi'(T')=\alpha_k\circ(\alpha_{k'-1}'
\circ\dots\circ\alpha_1')(T')$$
By lemma \ref{L1} we get from here that 
$$(\alpha_{k-1}\circ\dots\circ\alpha_1)(T)
=(\alpha_{k'-1}'
\circ\dots\circ\alpha_1')(T')$$
and the inductive hypothesis implies that 
$(\alpha_{k-1}'\circ\dots\circ\alpha_1') = 
(\alpha_{k-1}\circ\dots\circ\alpha_1)\circ\theta$
and $T=\theta(T')$ for 
$\theta=\alpha_{k'-k}'\circ\dots\circ\alpha_1'$
Therefore also $\varphi'=\varphi\circ\theta$. 

It remains to consider the case when $\alpha_k$ is a 
projection. By our assumption this is possible only if $k>1$.  
Moreover, $\alpha_{k-1}$ cannot be a projection. One can 
check directly  that $\alpha_{k'}'=\alpha_k$ and 
$\alpha_{k'-1}'=\alpha_{k-1}$, and  then argue 
the same way as above to complete the proof.  
\end{proof}

\begin{lemma}
\label{L3}
Let $T\in s\bdn^{k+1}\setminus\bdn^{k}$ and let 
$\varphi=\alpha_m\circ\dots\circ\alpha_1$ be a morphism 
in $\D$ such that $\alpha_1$ is not a projection. 
Then $\varphi(T)\not\in\bdn^k$.
\end{lemma} 

\begin{proof}
Notice that by the definition of $\bdn^k$ and $s\bdn^k$
if $S$ is an element of $\bdn^k$ then $S=\psi(S')$ for some 
$\psi\in\D$ and $S'\in s\bdn^k$. We can also assume that 
$\psi = \alpha_l'\circ\dots\circ\alpha_1'$ where 
$\alpha_1'$ is not a projection. 

Assume that for some $T\in s\bdn^{k+1}\setminus\bdn^{k}$ 
and for some $\varphi\in\D$
we have $\varphi(T)\in\bdn^k$. 
Then $\varphi(T)=\varphi'(T')$  for some $\varphi'\in\D$
and $T'\in s\bdn^k$. Thus by (\ref{L2})  
either $T'=\theta{T}$ or $T=\theta(T')$ for some 
morphism $\theta\in\D$. The latter is impossible:
indeed, we would get $T=\theta(T')\in\bdn^k$ 
contrary to the assumption we made about $T$. Therefore 
$T'=\theta(T)$ and we can also assume that the decomposition of 
$\theta$ into generators of $\D$ does not start with 
a projection. Replacing $\theta$ with $p\circ\theta$ 
if necessary (where $p$ is a projection morphism) we get
that $\theta(T)\in\bdn^{k-1}$. Continuing with this argument we 
would eventually have to conclude that there exists a morphism
$\theta'=\alpha_l'\circ\dots\circ\alpha_1'\in\D$ such that 
$\alpha_1'$ is not a projection morphism and $\theta'(T)$
belongs to $\bdn^0$. But $\bdn^0=D_n$, so $\theta'(T)$
would have to be represented by a morphism of $\D$. 
One can check that this can happen only if $T$ itself 
belongs to $D_n$. But $D_n\subseteq\bdn^k$, so this last 
statement contradicts our assumption on $T$.     
\end{proof}

\begin{lemma}
\label{L4}
Let $T\in\bdn^{k+1}\setminus s\bdn^{k+1}$. There exists 
a unique element $S\in s\bdn^{k+1}\setminus\bdn^k$ and 
a unique morphism $\varphi\in\D$ such that 
$\varphi(S)=T$. 
\end{lemma}

\begin{proof}
Assume that $T=\varphi(S)=\varphi'(S')$ for some 
$\varphi, \varphi'\in\D$ and some  
$S, S'\in s\bdn^{k+1}\setminus\bdn^k$.
Since decompositions of $\varphi$ and $\varphi'$
into generators of $\D$ cannot start with a projection,
by  (\ref{L2}) we can assume that $S=\theta(S')$
and $\varphi'=\varphi\circ\theta$
for some morphism $\theta\in\D$. It follows that 
$\theta(S')\in s\bdn^{k+1}$ and composing $\theta$
with some projection morphism $p$ we have
$p\circ\theta(S')\in\bdn^{k}$. Lemma \ref{L3}
implies that the decomposition of $p\circ\theta$ 
into generators of $\D$ must begin with a projection.
Therefore $\theta$ must be an identity morphism. 
It follows that $S=S'$ and $\varphi=\varphi'$.

\end{proof}

Finally we can turn to 

\begin{Proofa}
{\em {\hskip -3mm} of lemma \ref{PULLBACK LEMMA}.\ }
We will show that for any $\D$--diagram $X$
the square
$$\xymatrix{
\Hom_{\D}(\bdn^k,X)\ar[d] & \Hom_{\D}(\bdn^{k+1}, X)\ar[l]\ar[d] \\
\Hom_{\Prj}(\bdn^k, X)     & \Hom_{\Prj}(s\bdn^{k+1}, X)\ar[l] \\
}$$ 
is a pullback diagram of sets. 
Since limits of diagrams of simplicial sets can be calculated 
by taking the limit in each simplicial dimension separately, 
and since for  
$Y,X\in\spaces^{\D}$ the set of $k$-dimensional simplices
$\map_{\D}(Y,X)$ can be described as 
$\map_{\D}(Y,X)_m=\Hom_{\D}(Y,X^{\Delta[m]})$.
the statement of lemma \ref{PULLBACK LEMMA} will follow.

Let $P$ denote the set which is the pullback of the diagram 
$$\Hom_{\D}(\bdn^k,X)\lra\Hom_{\Prj}(\bdn^k, X)
\longleftarrow\Hom_{\Prj}(s\bdn^{k+1},X)$$  
The commutativity of the square diagram above implies that 
there exists a unique map 
$\iota\colon \Hom_{\D}(\bdn^{k+1}, X)\ra P$
such that the following diagram commutes: 
$$\xymatrix{
 & &\Hom_{\D}(\bdn^{k+1}, X)\ar[dl]!U|{\iota}\ar@/_/[dll]\ar@/^/[ddl]\\
\Hom_{\D}(\bdn^k,X)\ar[d] & P\ar[l]\ar[d] & \\
\Hom_{\Prj}(\bdn^k, X)     & \Hom_{\Prj}(s\bdn^{k+1}, X)\ar[l]  & \\
}$$
We want to construct the inverse of the map $\iota$.
Notice that the elements of $P$ are maps of $\Prj$--diagrams
$\epsilon\colon s\bdn^{k+1}\ra X$ such that 
$\epsilon|_{\bdn^k}\colon\bdn^k\ra X$ commutes with all 
morphisms of $\D$. We claim that such $\epsilon$ admits a unique 
extension to a map of $\D$--diagrams 
$\bar{\epsilon}\colon\bdn^{k+1}\ra X$. To see that 
take $T\in \bdn^{k+1}\setminus s\bdn^{k+1}$. By (\ref{L4})
there exists a unique element $S\in s\bdn^{k+1}$ and 
$\varphi\in\D$ such that $\varphi(S)=T$. Define 
$\bar{\epsilon}(T):= \varphi(\epsilon(S))$. Uniqueness of 
this extension is obvious, and using (\ref{L3}) and
(\ref{L4}) it is also not difficult to check that 
$\bar{\epsilon}$ is a well defined  map of $\D$--diagrams.
The correspondence $\epsilon\mapsto\bar{\epsilon}$ gives a 
map $\nu\colon P\lra \Hom_{\D}(\bdn^{k+1}, X)$.
One can verify that the compositions $\nu\circ\iota$
and $\iota\circ\nu$ are both identities. 
Therefore we get $P\cong \Hom_{\D}(\bdn^{k+1}, X)$.
\end{Proofa}

%
%

\section{The pointed case}
\label{POINTED CASE}
We describe a few changes one needs to make 
in order to prove theorems  \ref{MAIN 1} and \ref{MAIN 2}
for pointed semi-theories. Since, as we noted in 
\S\ref{INTRO}, the arguments we used for unpointed semi-theories
apply in this case with minor changes only, we will concentrate
on the differences in some definitions and constructions.

First, since free categories usually do not have terminal
and initial objects, in the pointed case definition 
~\ref{DEF FREE SEMI} must be modified as follows:

\begin{definition}
A pointed semi-theory $\C$ is free if

\noindent -- there exists a free subcategory $\C_{>0}\subset\C$
with objects $[1],[2],\dots$ such that all projections 
in $\C$ are free generators of $\C_{>0}$;

\noindent -- the only morphisms of $\C$ which do not belong to 
$\C_{>0}$ are the morphims $[n]\ra [0]$, $[0]\ra [n]$ for 
$n\geq 0$ and their compositions.
 
\end{definition}

\noindent In this setting by generators of a free 
pointed semi-theory $\C$ we understand the generators 
of $\C_{>0}$.

Notice that the category $\C_{>0}$ is a free unpointed 
semi-theory. Going in the opposite direction,  
given any free unpointed semi-theory $\C$ one can construct 
a unique free pointed semi-theory $\C_+$ such that 
$(\C_+)_{>0}=\C$. Applying this construction 
to the semi-theory $\Prj$ (sec. \ref{INITIAL}) we 
obtain a semi-theory $\Prj_+$ which is initial 
among all pointed semi-theories.   

For a free pointed semi-theory $\C$
the algebraic theory $\CC$ is constructed in a similar 
way as in the unpointed case. A tree $T$ 
representing a morphisms $[n]\ra[1]$ in $\CC$ has its non-initial  
edges labeled with generators of $\C$. However, 
one must allow that the labels of the initial edges can be 
the projection morphisms and the morphism $\iota_n\in\C$ 
which is the (unique) composition 
$[n]\ra [0] \ra [1]$. As a consequence morphisms in 
$\Hom_\CC([n],[1])$ are now equivalence classes of trees 
rather than trees themselves. The equivalence relation 
is generated by the following condition: trees  
$T$ and $T'$ are equivalent if $T'$ is obtained by removing
an initial edge of $T$ labeled with $\iota_n$ and grafting
in its place any tree $S$ whose all initial edges are labeled 
with $\iota_n$. 

The composition of $(T_1,\dots, T_m)$ 
representing a morphisms in $\Hom_\CC([n],[m])$ with 
$S\in\Hom_\CC([m],[1])$ is a morphisms represented 
by the tree obtained by grafting $T_k$ in place of each
edge of $S$ labeled with $p^m_k$ and relabeling with $\iota_n$
all initial edges of $S$ which have the label $\iota_m$.

Finally, the morphisms 
in $\CC$ with the source or target $[0]$ are determined by 
the conditions that $[0]$ is the initial and the terminal
object in $\CC$, and that the composition $[n]\ra [0]\ra [1]$
in $\CC$ is represented by the tree which has only one edge, 
and it is labeled with $\iota_n$. 

In order to construct the resolution $\FC$ for a pointed
semi-theory $\C$ take $\bo F_0\C$ to be the free 
pointed semi-theory such that $\bo (F_0 \C)_{>0}$ is the free
category generated by the full subcategory of $\C$ on 
objects $[1],[2], \dots$. For $k>0$ we construct
$\bo F_{k+1} \C$ out of $\bo F_k \C$ in the same manner. 
The algebraic theory $\bFC$ is obtained by 
applying the (pointed) completion in every simplicial 
dimension of $\FC$.

%
%

\end{document}